\documentclass[11pt,leqno]{article}

\usepackage[english]{babel}

\usepackage{amsmath}

\usepackage{amssymb}

\usepackage{xypic}

\usepackage[all]{xy}

\usepackage{latexsym}

\usepackage{amsfonts}

\usepackage{epsf}

\usepackage[dvips]{graphicx}

\usepackage[latin1]{inputenc}

\usepackage[T1]{fontenc}

\numberwithin{equation}{section}

\newcommand{\T}{\mathcal{T}}

\renewcommand{\mod}{\mathrm{Mod}}

\newcommand{\R}{\mathbb{R}}

\newcommand{\Q}{\mathbb{Q}}

\newcommand{\N}{\mathbb{N}}

\newcommand{\I}{\mathrm{I}}

\newcommand{\iso}{\stackrel{\sim}{\to}}

\newcommand{\rh}{\mathit{R}\mathcal{H}\mathit{om}}

\newcommand{\ho}{\mathcal{H}\mathit{om}}

\newcommand{\Ho}{\mathrm{Hom}}

\newcommand{\Rh}{\mathrm{RHom}}

\renewcommand{\dim}{\textbf{Proof.}}

\newcommand{\qed}{\nopagebreak \phantom{} \hfill $\Box$ \\}

\newcommand{\imin}[1]{#1^{-1}}

\newcommand{\Lpro}{\underleftarrow{\lim}}

\newcommand{\lpro}[1]{\underset{#1}{\underleftarrow{\lim}}}

\newcommand{\exs}[3]{0 \to {#1} \to {#2} \to {#3} \to 0}

\newtheorem{teo}{Theorem}[section]

\newtheorem{df}[teo]{Definition}

\newtheorem{oss}[teo]{Remark}

\newtheorem{prop}[teo]{Proposition}

\title{\bf{ON THE HOMOLOGICAL DIMENSION OF O-MINIMAL AND SUBANALYTIC SHEAVES}}
\date{}

\usepackage{latexsym}
\begin{document}

\maketitle

\begin{abstract}
Here we prove that the homological dimension of the category of sheaves on a topological space satisfying some suitable conditions is finite. In particular, we find conditions to bound the homological dimension of o-minimal and subanalytic sheaves.
\end{abstract}

\section{Introduction}

In \cite{lucap} we studied the category $\mod(k_{X_{sa}})$ of sheaves on the subanalytic site $X_{sa}$ associated to a real analytic manifold $X$. We defined the subcategory of quasi-injective sheaves (i.e. $F$ is quasi-injective if the restriction $\Gamma(U;F)\to \Gamma(V;F)$ is surjective for each $U \supseteq V$ open subanalytic relatively compact) and we saw that quasi-injective are injective with respect to the functors of direct image, proper direct image and $\ho(G,\cdot)$ when $G$ is a $\R$-constructible sheaf on $X$. Moreover we proved that the quasi-injective dimension of $\mod(k_{X_{sa}})$ is finite, which implies that the cohomological dimension of the above functors is finite. However we had no answer concerning the homological dimension of $\mod(k_{X_{sa}})$ (see Remark 2.3.5 of \cite{lucap}). The aim of this paper is to show that under some conditions (concerning the cardinality $2^{\aleph_0}$ of the reals) it is possible to bound this dimension. The key point is the fact that locally the category of subanalytic sheaves is equivalent to the category of sheaves on a topological space $\widetilde{X}$. We are reduced to work with $k$-sheaves on a topological space, and in this case the homological dimension is equal to the flabby dimension. Hence we are reduced to bound the cohomological dimension of $\Gamma(U;\cdot)$, for any open subset $U$ of $\widetilde{X}$. In order to do that we need to assume that the cardinality of subanalytic subsets (which is equal to $2^{\aleph_0}$) is smaller than $\aleph_k$, $k<\infty$.\\

In more details the contents of this paper are as follows.

In $\S$\,\ref{0} we recall some notions as the definitions of injective and flabby sheaf and the homological and flabby dimension of the category of sheaves. We refer to \cite{ks1} for more details.

In $\S$\,\ref{1} we study the general case of sheaves on a topological space $X$ with a basis $\T$ whose elements are Lindel\"of and stable under finite unions and intersections. We define the subcategory of $\T$-flabby sheaves on $X$ by saying that $F$ is $\T$-flabby if the restriction $\Gamma(U;F)\to \Gamma(V;F)$ is surjective for each $U,V\in\T$ with $U \supseteq V$. Since $\T$ forms a basis of the topology of $X$, we have that each open $U$ of $X$ has a covering $\{U_i\}_{i\in I}$ with $U_i\in\T$. Hence
$$
R\Gamma(U;F) \simeq R\lpro {i\in I}R\Gamma(U_i;F),
$$
and to bound the cohomological dimension of the functor $\Gamma(U;\cdot)$ is sufficient to bound the cohomological dimension of the functors $\Gamma(V;\cdot)$, $V \in \T$ and the cohomological dimension of the projective limit. The cohomological dimension of the projective limit is bounded if the cardinality of the index set $I$ is smaller than $\aleph_k$ for some $k<\infty$, (see \cite{jen} or \cite{fp} in the more general setting of quasi-abelian categories). Hence it is bounded if the cardinality of $\T$ is smaller than $\aleph_k$, $k<\infty$.

 In $\S$\,\ref{2} we consider an o-minimal structure ${\cal M}=(M,<,\ldots)$. Let $\widetilde{X}$ be the o-minimal spectrum of a definable space $X$. In this case $\T$ is the family of open constructible subsets of $\widetilde{X}$. In the case of an o-minimal expansion of an ordered group, by a result of \cite{ejp}, the $\T$-flabby dimension of $\widetilde{X}$ is finite. Moreover the cardinality of $\T$ is bounded by the product of the cardinality of $M$ and the cardinality of the language of the structure ${\cal M}$.

The case of subanalytic sheaves on a real analytic manifold is studied in $\S$\,\ref{3}. We first reduce to $\mod(k_{U_{X_{sa}}})$ where $U$ is a relatively compact subanalytic subset of $X$ isomorphic to $\R^N$ endowed with the Grothendieck topology induced by $X$. This category is equivalent to $\mod(k_{\widetilde{U}})$, where $\widetilde{U}$ is the o-minimal spectrum of $U$. In this case $\T$ is the family of open globally subanalytic subsets of $U$, and its cardinality is $2^{\aleph_0}$. Hence if we assume that $2^{\aleph_0}$ is smaller or equal than $\aleph_k$, $k<\infty$ (in the case $k=1$ this is the continuum hypothesis), then the homological dimension of $\mod(k_{\widetilde{U}})$ (and hence of $\mod(k_{U_{X_{sa}}})$) is bounded.

We end this work with an example showing that in $\widetilde{X}$ there are open subsets which do not admit countable covers.

\section{Some preliminaries}\label{0}

We introduce some fundamental results about sheaves we will use in the rest of the paper. We refer to \cite{ks1} for more details. Let $X$ be a topological space and let $k$ be a field. As usual, we denote by $\mod(k_X)$ the category of sheaves of $k$-vector spaces.

\begin{df} Let $F \in \mod(k_X)$.
\begin{itemize}
\item[-] $F$ is injective if the functor $\Ho(\cdot,F)$ is exact on $\mod(k_X)$.
\item[-] $F$ is flabby if for any open subset $U$ of $X$ the restriction morphism $\Gamma(X;F) \to \Gamma(U;F)$ is surjective.
\end{itemize}
\end{df}

In general injective $\Rightarrow$ flabby. When $k$ is a field we have injective $\Leftrightarrow$ flabby (\cite{ks1}, Exercise II.10).

\begin{df} The homological (resp. flabby) dimension of the category $\mod(k_X)$ is
the smallest $N \in \N \cup \{\infty\}$ such that for any $F \in
\mod(k_X)$ there exists an exact sequence
$$0 \to F \to \I^0 \to \cdots \to I^N \to 0$$
with $I^j$ injective (resp. flabby) for $0 \leq j \leq N$.
\end{df}

We shall need the following results.
\begin{itemize}
\item[-] The homological dimension of $\mod(k_X)$ is finite if and only if there exists $N \in \N$ such that $R^j\Ho(G,F)=0$ for any $F,G \in \mod(k_X)$ and any $j > N$.
\item[-] The flabby dimension of $\mod(k_X)$ is finite if and only if there exists $N \in \N$ such that for any open subsets $U$ of $X$, any $F \in \mod(k_X)$ and any $j > N$ we have $R^j\Gamma(U;F)=0$ (i.e. if the functor $\Gamma(U;\cdot)$ has finite cohomological dimension).
\end{itemize}

In particular, when  $k$ is a field, since injective $\Leftrightarrow$ flabby, the homological dimension is finite if and only if the functor $\Gamma(U;\cdot)$ has finite cohomological dimension for any open subsets $U$ of $X$.

\section{$\T$-flabby sheaves and homological dimension}\label{1}

\begin{df} A Lindelöf space is a topological space in which every open cover has a countable subcover.
\end{df}

\begin{df} The Lindelöf degree $l(X)$ of a topological space $X$ is the smallest cardinal $\kappa$ such that every open cover of the space X has a subcover of size at most $\kappa$. In this notation, X is Lindelöf iff $l(X) = \aleph_0$.
\end{df}

Let $k$ be a field. Let $X$ be a topological space and suppose that $X$ admits a family $\T$ of open subsets such that
\begin{itemize}
\item[T1:] each element of $\T$ is Lindel\"of,
\item[T2:] $\T$ is stable under finite unions and intersections,
\item[T3:] $\T$ forms a basis for the topology of $X$.
\end{itemize}

\begin{df} A sheaf $F \in \mod(k_X)$ is said to be $\T$-flabby if the restriction morphism $\Gamma(U;F) \to \Gamma(V;F)$ is surjective for each $U,V \in \T$ with $V \subseteq U$.
\end{df}

Remark that flabby $\Rightarrow$ $\T$-flabby.

\begin{df} The $\T$-flabby dimension of $\mod(k_X)$ is the smallest $N \in \N \cup \infty$
such that for any $F \in
\mod(k_X)$ there exists an exact sequence
$$0 \to F \to \I^0 \to \cdots \to I^N \to 0$$
with $I^j$ $\T$-flabby for $0 \leq j \leq N$.
\end{df}

\begin{prop}
\label{prop:soft_section_exact}
  Let $\exs{F'}{F}{F''}$ be an exact
  sequence in $\mod(k_X)$ with $F'$ $\T$-flabby. Then for any open
  subset $U$ which is Lindel\"of the morphism
$$
\Gamma(U;F) \to \Gamma(U;F'')
$$
is
surjective.
\end{prop}
\dim\ \
  We first consider a section $s''\in \Gamma(U;F'')$.  Since $F \to F''$ is surjective, we may find a covering $U = \bigcup_{i\in\I} U_i$, $U_i \in \T$ and $s_i \in
  \Gamma(U_i;F)$ whose image is $s''|_{U_i}$. Since $U$ is Lindel\"of, we may find a countable subcover
  $U = \bigcup_{n\in\N} U_n$, $U_n \in \T$ and $s_n \in
  \Gamma(U_n;F)$ whose image is $s''|_{U_n}$.

  Set $V_n = \bigcup_{i=1}^nU_i$. We prove by induction
  on $n$ that there exists a section $t_{n+1} \in \Gamma(V_{n+1};F)$
  whose image is $s''|_{V_{n+1}}$ and $t_{n+1} |_{V_n}= t_n$.

  This is clear for $n=0,1$ and we assume it is proved for $n$. By the induction hypothesis there exists a section $t_n \in \Gamma(V_n;F)$
  whose image is $s''|_{V_n}$ and $t_n |_{V_m}
  = t_m$ if $m<n$. We set for short $V_n=W_1$ and $U_{n+1}=W_2$.
  We have seen that there exist $t_j \in
\Gamma(W_j;F)$ whose image is $s''|_{W_j}$ for $j=1,2$. On $W_1
\cap W_2$ $t_1-t_2$ defines a section of $\Gamma(W_1 \cap W_2;F')$
which extends to $t' \in \Gamma(V_{n+1};F')$ because $F'$ is $\T$-flabby. Replace $t_2$ with
$t_2+t'$. We may suppose that $t_1=t_2$ on $W_1 \cap W_2$. Then
there exists $t_{n+1} \in \Gamma(V_{n+1};F)$ such that
$t_{n+1}|_{W_j}=t_j$, $j=1,2$. The
  $t_n$'s glue together into a section
  $$
  s \in \lpro {n\in\N}\Gamma(V_n;F) \simeq \Gamma(U;F)
  $$
  which is sent to $s''$, which proves the surjectivity of
  the morphism.
\qed

\begin{prop}\label{www}
Let $F',F,F'' \in \mod(k_X)$, and consider the
exact sequence
$$
\exs{F'}{F}{F''}.
$$
Suppose that $F'$ is $\T$-flabby. Then $F$ is $\T$-flabby if and only if $F''$ is $\T$-flabby.
\end{prop}
\dim\ \ Let $U,V\in \T$  with $V \subseteq U$ and let us
consider the diagram below
$$ \xymatrix{0 \ar[r] & \Gamma(U;F') \ar[d]^\alpha \ar[r] & \Gamma(U;F) \ar[d]^\beta \ar[r] & \Gamma(U;F'') \ar[d]^\gamma \ar[r] & 0\\
0 \ar[r] & \Gamma(V;F') \ar[r] & \Gamma(V;F) \ar[r] & \Gamma(V;F'') \ar[r] & 0} $$
where the row are exact by Proposition \ref{prop:soft_section_exact} and
the morphism
$\alpha$ is surjective since $F'$ is $\T$-flabby. It follows from the five lemma that $\beta$ is surjective if and only if $\gamma$ is surjective.
\qed

\begin{prop}\label{TUinj} Let $F$ be $\T$-flabby. Then $F$ is $\Gamma(U;\cdot)$-injective for each open $U\subseteq X$ which is Lindel\"of.
\end{prop}
\dim\ \ The family of $\T$-flabby sheaves contains injective sheaves, hence it is cogenerating. Then the result follows from Propositions \ref{prop:soft_section_exact} and \ref{www}.
\qed

Now let us consider the case of open subsets of Lindel\"of degree $\aleph_k$, $k<\infty$. In order to do that we need the following result.

\begin{prop}\label{lprodim} Let $X$ be a topological space admitting a family $\T$ satisfying T1-T3. Suppose that $X$ has finite $\T$-flabby dimension. Let $U$ be an open subset of $X$ of Lindel\"of degree $\aleph_k$, $k<\infty$. Then the cohomological dimension of $\Gamma(U;\cdot)$ is finite.
\end{prop}
\dim\ \ We have
$$
R\Gamma(U ; F) \simeq R\lpro {i \in I}R\Gamma(U_i; F)
$$
where $U=\bigcup_{i\in I}U_i$ with $U_i \in \T$ and $\sharp I$ smaller than $\aleph_k$ . Since $X$ has finite $\T$-flabby dimension and using Proposition \ref{TUinj} for each $i$ we may replace $R\Gamma(U_i;F)$ with $\Gamma(U_i;I^\bullet)$, where $I^\bullet$ is a $\T$-flabby resolution of $F$ of length $N<\infty$. Since
cohomological dimension of $\lpro {i \in I}$
is finite if $\sharp I$ is smaller than $\aleph_k$, $k<\infty$ (see \cite{jen,fp}), then the $j$-th cohomology of $R\lpro {i \in I}R\Gamma(U_i; F)$ vanishes for $j>N+M$. Since $M,N<\infty$ are independent of $F$ and $i$, the result follows.
\qed

\begin{teo} \label{homdim} Let $X$ be a topological space admitting a family $\T$ satisfying T1-T3. Suppose that $X$ has finite $\T$-flabby dimension. Suppose that there exists $k<\infty$ such that every open subset of $X$ has Lindel\"of degree $\aleph_k$. Then $\mod(k_X)$ has finite homological dimension.
\end{teo}
\dim\ \ By Proposition \ref{lprodim} $X$ has finite flabby dimension. In the case of sheaves of $k$-vector spaces we have flabby $\Leftrightarrow$ injective and the result follows.
\qed

\section{Homological dimension of o-minimal sheaves}\label{2}

Let ${\cal M}=(M,<,\ldots)$ be an o-minimal structure. In \cite{ejp} the authors studied sheaf cohomology of sheaves on $\widetilde{X}$, the o-minimal spectrum of a definable space $X$. The category of sheaves on $\widetilde{X}$ is equivalent to the category of sheaves on the o-minimal site, consisting of open definable subsets of $X$ and coverings admitting a finite refinement.

Let $\T$ be the family of open constructible subsets of $\widetilde{X}$. Then $\T$ satisfies T1-T3, indeed
\begin{itemize}
\item[-] $\T$ forms a basis for the topology of $\widetilde{X}$,
\item[-] every element of $\T$ is quasi-compact,
\item[-] $\T$ is stable under finite unions and intersections.
\end{itemize}
For each open subset $U$ of $\widetilde{X}$ and each $F \in \mod(k_{\widetilde{X}})$ we have
$$
R\Gamma(U;F) \simeq R\lpro {i\in I}R\Gamma(U_i;F)
$$
where $U=\bigcup_{i \in I}U_i$ and $U_i \in \T$ for each $i \in I$. In order to apply Theorem \ref{homdim} we need that:
\begin{itemize}
\item[-] the $\T$-flabby dimension of $\widetilde{X}$ is finite,
\item[-] $\sharp I \leq \aleph_k$, $k<\infty$.
\end{itemize}
If every open constructible subset $U$ of $X$ is normal (this is true in the case of an o-minimal expansion of an ordered group), then by Proposition 4.2 of \cite{ejp} we have $R^j\Gamma(U;F)=0$ for all $j>{\rm dim}U={\rm dim}X$.

Concerning $\sharp I$, we have that $\sharp I \leq \sharp \T$, and $\sharp \T \leq \sharp {\cal L} \cup M=\sharp {\cal L} \cdot \sharp M$, where ${\cal L}$ denotes the language of the structure ${\cal M}$.

\section{Homological dimension of subanalytic sheaves}\label{3}

Let $X$ be a real analytic manifold, let $X_{sa}$ the associated subanalytic site and denote by $\mod(k_{X_{sa}})$ the category of sheaves of $k$-vector spaces on $X_{sa}$. We refer to \cite{ks2,lucap} for the theory of subanalytic sheaves.
We will see that under suitable hypothesis the homological dimension of $\mod(k_{X_{sa}})$ is finite, i.e. there exists $N \in \N$ such that $R^j\Ho(G,F)=0$ for any $F,G \in \mod(k_{X_{sa}})$ and any $j > N$.


Let $X$ be a real analytic manifold. Then $X$ has an atlas $(U_n,\phi_n)_{n\in\N}$, $\phi_n:U_n \iso \R^M$, $U_n$ relatively compact open subanalytic subset.
Let us consider the subanalytic site $X_{sa}$ associated to $X$. In order to bound $\Rh(F,G)$ for any $F,G \in \mod(k_{X_{sa}})$ it suffices to bound $\rh(F,G)$. Indeed we have $\Rh(F,G) \simeq R\Gamma(X;\rh(F,G))$ and $\Gamma(X;\cdot)$ has finite cohomological dimension (see \cite{lucap}).

We have $R^j\ho(F,G)=0$ if $R^j\ho(F,G)|_{U_n}=0$ for each $n\in\N$. Then we may reduce to the case $U_n$ with the Grothendieck topology induced by the subanalytic site, i.e. the category of globally subanalytic open subsets of $\R^M$ with coverings admitting a finite subcover.
Globally subanalytic geometry is defined by the o-minimal structure given by the ordered field of real numbers expanded by globally analytic functions  ${\mathbb R}_{{\rm an}}=({\mathbb R}, <, 0,1, +, \cdot, (f)_{f\in {\rm an}})$, where $f$ is a restriction to $[-1,1]^m$ of a convergent power series on some neighborhood of $[-1,1]^m$.


Set $X=U_n$ and let $\widetilde{X}$ be the o-minimal spectrum of $X$.
The family $\T$ of open sets $\widetilde{U}$ where $U$ is
open globally subanalytic satisfies T1-T3. Moreover the $\T$-flabby dimension of $\mod(k_{\widetilde{X}} )$ is finite.
Indeed by Proposition 4.2 of \cite{ejp} the cohomological dimension of $\Gamma(\widetilde{U};\cdot)$ is finite for each $\widetilde{U} \in \T$.

Let us see that $\sharp\T \leq 2^{\aleph_0}$. In order to see that we have to show that $\sharp {\cal L} = 2^{\aleph_0}$, where ${\cal L}$ denotes the language of $\R_{{\rm an}}$. This is true if $\sharp\{(f)_{f\in{\rm an}}\}=2^{\aleph_0}$. This is seen by identifying $f=\sum_{I \in \N^m}a_I x^I$ with $(a_I)_{I \in \N^m} \in \R^{\N^m}$. Then $\sharp \{(f)_{f\in{\rm an}}\} \leq \sharp \R^{\N^m} = (2^{\aleph_0})^{\aleph_0} = 2^{\aleph_0 \cdot \aleph_0} = 2^{\aleph_0}$ (see \cite{jech}).

Let $U$ be an open subset of $\widetilde{X}$. Then $U =\bigcup_{i\in I} U_i$ with $U_i \in \T$
and $\sharp I$ is smaller than $\sharp \T = 2^{\aleph_0}$. Hence every open subset of $X$ has Lindel\"of degree at most $2^{\aleph_0}$. If we assume that $2^{\aleph_0}$ is smaller than
$\aleph_k$ for some $k<\infty$ (if we suppose $\leq \aleph_1$ it is the continuum hypothesis), then Theorem \ref{homdim} implies that the homological dimension of $\mod(k_{\widetilde{X}} )$ 
is finite.

\begin{oss}\label{lindim}
In $\widetilde{X}$ there are open subsets which are not Lindel\"of. For example let us consider the open set (in $\widetilde{X}=\widetilde{\R^2}$ with respect to $\R_{\mathrm{an}}$)
$$
U=\bigcup_{r \in (0,1) \setminus \Q}\widetilde{V}_r,
$$
where
$$
V_r=\{(x,y) \in \R^2\;;\;0<x<1,\;r<y<r+x\}.
$$
We prove that $\{V_r\}_{r \in (0,1)\setminus\Q}$ has no countable subcover.\\

We argue by contradiction.  Suppose that there exists a countable subcover $\{\widetilde{V}_{r_n}\}_{n\in\N}$ of $U$ in $\widetilde{X}$. Let $0^+$ be the ultrafilter defined by
$0^+=\{S \supset (0,\varepsilon)\}$,
where $\varepsilon>0$ and $S$ is globally subanalytic. Let $\pi:\R^2 \to \R$ be the projection onto the first coordinate.
We have
$$
U \cap \imin {\widetilde{\pi}}(0^+) = \bigsqcup_{r\in (0,1)\setminus \Q} \left( \widetilde{V}_r \cap \imin {\widetilde{\pi}}(0^+) \right).
$$
Indeed, let $x \in \widetilde{V}_r \cap \imin {\widetilde{\pi}}(0^+)$, and let $s \neq r$.
Let $\varepsilon < |r-s|$, then $V_s \cap \imin \pi((0,\varepsilon)) \cap V_r = \emptyset$ and $\widetilde{V_s} \cap \imin {\widetilde{\pi}}(0^+) \subset \widetilde{V_s} \cap \imin {\widetilde{\pi}}(\widetilde{(0,\varepsilon)})$ since $0^+ \subset \widetilde{(0,\varepsilon)}$ for any $\varepsilon>0$. Hence $\widetilde{V}_r \cap \imin {\widetilde{\pi}}(0^+) \cap \widetilde{V}_s = \emptyset$ if $r \neq s$.
Moreover
$$
\bigsqcup_{r\in (0,1)\setminus \Q} \left( \widetilde{V}_r \cap \imin {\widetilde{\pi}}(0^+) \right) \subset \bigcup_{n\in\N}\widetilde{V}_{r_n}=U.
$$
Since $\sharp((0,1)\setminus\Q)=2^{\aleph_0}>\aleph_0$, there exists $n\in\N$ such that
$$
\sharp\left\{t \in (0,1)\setminus\Q \;;\; \widetilde{V_t} \cap \imin {\widetilde{\pi}}(0^+) \subseteq \widetilde{V}_{r_n} \right\} \geq 2.
$$
But $\widetilde{V_t} \cap \imin {\widetilde{\pi}}(0^+) \subseteq \widetilde{V}_{r_n}$ implies that $V_t=V_{r_n}$, so $t=r_n$ which is a contradiction.
\end{oss}

\begin{oss} 
More generally, as in Remark \ref{lindim}, given an o-minimal expansion of an ordered group $(\R,<,+,\ldots)$, for any definable manifold $X$ of dimension $\geq 2$, the space $\widetilde{X}$ has open subsets which have Lindel\"of degree $2^{\aleph_0}$.
\end{oss}

\begin{oss} It seems that the right setting is that of locally definable
spaces in an o-minimal structure. In that setting the o-minimal (definable)
case and the subanalytic case would be treated uniformly. However, to work
on that setting one has to take into account the (special case of the) theory
locally semi-algebraic spaces developed in \cite{D3,dk5}.
\end{oss}

\end{document}